\input amstex
\input amsppt.sty
\magnification=\magstep1
\hsize=32truecc
\vsize=22.2truecm
\baselineskip=16truept
\NoBlackBoxes
\TagsOnRight \pageno=1 \nologo

\def\l{\left}
\def\r{\right}
\def\bg{\bigg}
\def\({\bg(}
\def\[{\bg\lfloor}
\def\){\bg)}
\def\]{\bg\rfloor}
\def\t{\text}
\def\f{\frac}

\def\p{\ (\roman{mod}\ p)}

\def\bi{\binom}
\def\eq{\equiv}

\def\ls{\leqslant}

\def\mo{\roman{mod}}

\def\Proof{\noindent{\it Proof}}

\def\Remark{\medskip\noindent{\it  Remark}}

\def\Ack{\medskip\noindent {\bf Acknowledgment}}
\hbox {J. Number Theory 134(2014), no.\,1, 181--196.}
\bigskip
\topmatter
\title $p$-adic congruences motivated by series\endtitle
\author Zhi-Wei Sun\endauthor
\leftheadtext{Zhi-Wei Sun} \rightheadtext{$p$-adic congruences motivated by series}
\affil Department of Mathematics, Nanjing University\\
 Nanjing 210093, People's Republic of China
  \\  zwsun\@nju.edu.cn
  \\ {\tt http://math.nju.edu.cn/$\sim$zwsun}
\endaffil
\abstract Let $p>5$ be a prime. Motivated by the known formulae
$$\sum_{k=1}^\infty\f{(-1)^k}{k^3\bi{2k}k}=-\f 25\zeta(3)\ \ \t{and}\ \ \sum_{k=0}^\infty\f{\bi{2k}k^2}{(2k+1)16^k}=\f{4G}{\pi}$$
(where $G=\sum_{k=0}^\infty(-1)^k/(2k+1)^2$ is the Catalan constant), we show that
$$\gather\sum_{k=1}^{(p-1)/2}\f{(-1)^k}{k^3\bi{2k}k}\eq-2B_{p-3}\pmod p,
\\\sum_{k=(p+1)/2}^{p-1}\f{\bi{2k}k^2}{(2k+1)16^k}\eq-\f 74p^2B_{p-3}\pmod{p^3}
\endgather$$
and
$$\sum_{k=0}^{(p-3)/2}\f{\bi{2k}k^2}{(2k+1)16^k}\eq-2q_p(2)-p\,q_p(2)^2+\f 5{12}p^2B_{p-3}\pmod{p^3},$$
where $B_0,B_1,B_2,\ldots$ are Bernoulli numbers and $q_p(2)$ is the Fermat quotient $(2^{p-1}-1)/p$.
\endabstract
\thanks 2010 {\it Mathematics Subject Classification}.\,Primary 11B65;
Secondary 05A10, 05A19, 11A07, 11B68, 11S99.
\newline\indent Supported by the National Natural Science
Foundation (grant 11171140) of China.
\endthanks
\endtopmatter
\document

\heading{1. Introduction}\endheading

Let $p>3$ be a prime.
In 2010 the author and R. Tauraso [ST] proved that
$$\sum_{k=1}^{p-1}\f{\bi{2k}k}k\eq\f 89p^2B_{p-3}\pmod{p^3},$$
where $B_0,B_1,B_2,\ldots$ are Bernoulli numbers.
Note that
$$\bi{2k}k=\f{(2k)!}{k!^2}\eq0\pmod{p}\quad\t{for}\ k=\f{p+1}2,\ldots,p-1.$$
The author [S11c] managed to show that
$$\sum_{k=1}^{(p-1)/2}\f{\bi{2k}k}k\eq(-1)^{(p+1)/2}\f 83pE_{p-3}\pmod{p^2}$$
and
$$\sum_{k=1}^{(p-1)/2}\f1{k^2\bi{2k}k}\eq(-1)^{(p-1)/2}\f 43E_{p-3}\pmod{p},$$
where $E_0,E_1,E_2,\ldots$ are Euler numbers. It is interesting to compare the last congruence with the known formula
$$\sum_{k=1}^\infty\f1{k^2\bi{2k}k}=\f{\zeta(2)}3=\f{\pi^2}{18}.$$
Note that van Hamme [vH] and his followers ever considered $p$-adic analogues
of some hypergeometric series related to the Gamma function or $\pi=\Gamma(1/2)^2$
but Bernoulli numbers or Euler numbers never appeared in their work.

In 1979 Ap\'ery (cf. [Ap] and [vP]) proved the irrationality of
$\zeta(3)=\sum_{k=1}^\infty1/n^3$ and his following formula
$$\sum_{k=1}^\infty\f{(-1)^k}{k^3\bi{2k}k}=-\f25\zeta(3)$$
plays an important role in the proof. Motivated by this, Tauraso
[T10] showed that if $p>5$ is a prime then
$$\sum_{k=1}^{p-1}\f{(-1)^k}{k^3\bi{2k}k}\eq-\f 25\cdot\f{H_{p-1}}{p^2}\ (\mo\ p^3)$$
and $$\sum_{k=1}^{p-1}\f{(-1)^k}{k^2}\bi{2k}k\eq\f 45\cdot\f{H_{p-1}}p\ (\mo\ p^3),$$
where $H_n=\sum_{0<k\ls n}1/k\ (n=0,1,2,\ldots)$ are harmonic
numbers. It is well known (cf. [G1] or [Su, Theorem 5.1(a)]) that
$$H_{p-1}\eq-\f{p^2}3B_{p-3}\pmod{p^3}\quad \ \t{for any prime}\ p>3.$$
Actually Tauraso obtained $\sum_{k=1}^{p-1}\f{(-1)^k}{k^2}\bi{2k}k$ mod $p^4$ for each prime $p>5$ via putting $n=p$ in the following identity
$$\sum_{k=1}^n\bi{2k}k\f{k^2}{4n^4+k^4}\prod_{j=1}^{k-1}\f{n^4-j^4}{4n^4+j^4}=\f{2}{5n^2}$$
conjectured by J. M. Borwein and D. M. Bradley [BB] and proved by G. Almkvist and
A. Granville [AG].

Let $p>3$ be a prime. The author [S11c] proved that
$$\sum_{k=1}^{p-1}\f{\bi{2k}k^2}{16^k}\eq (-1)^{(p-1)/2}-p^2E_{p-3}\pmod{p^3}.$$
Recently Tauraso [T12] showed that
$$\sum_{k=1}^{p-1}\f{\bi{2k}k^2}{k16^k}\eq-2H_{(p-1)/2}\ (\mo\ p^3).$$
 Inputting the  command
\newline\indent \quad {\tt FullSimplify[Sum[Binomial[2k,k$]^\wedge2$/(k*16${}^\wedge$k),$\{$k,1,Infty$\}$]]},
\newline we get from {\tt Mathematica} (version 7) the identity
$$\sum_{k=1}^\infty\f{\bi{2k}k^2}{k16^k}=4\log2-\f{8G}{\pi},$$
where $G$ is the Catalan constant given by
$$G:=\sum_{k=0}^\infty\f{(-1)^k}{(2k+1)^2}.$$

Now we state our first theorem.

\proclaim{Theorem 1.1} Let $p>5$ be a prime. Then
$$\align\sum_{k=1}^{(p-1)/2}\f{(-1)^k}{k^3\bi{2k}k}\eq&-2B_{p-3}\pmod p,\tag1.1
\\\sum_{k=1}^{(p-1)/2}\f{(-1)^k}{k^2}\bi{2k}k\eq&\f{56}{15}pB_{p-3}\pmod{p^2},\tag1.2
\\\sum_{k=1}^{(p-1)/2}\f{\bi{2k}k^2}{k16^k}\eq&-2H_{(p-1)/2}-\f 72p^2B_{p-3}\pmod{p^3},\tag1.3
\\-\f 4{p^2}\sum_{k=(p+1)/2}^{p-1}\f{\bi{2k}k^2}{k16^k}\eq&\sum_{k=1}^{(p-1)/2}\f{16^k}{k^3\bi{2k}k^2}\eq -14B_{p-3}\pmod{p}.\tag1.4
\endalign$$
\endproclaim
\Remark\ 1.1. Both (1.1) and (1.2) were conjectured in [S11c, Conjecture 1.1]. The reader may consult [S10] and [S11a] for
$\sum_{k=0}^{p-1}\bi{2k}k/m^k$ modulo powers of $p$, where $p$ is an odd prime and $m$ is an integer not divisible by $p$.
\medskip

Motivated by the formulae
$$\sum_{k=0}^\infty\f{\bi{2k}k}{(2k+1)16^k}=\f{\pi}3\ \ \t{and}\ \ \sum_{k=0}^\infty\f{\bi{2k}k}{(2k+1)^2(-16)^k}=\f{\pi^2}{10},$$
the author [S11b] showed that
$$\sum_{k=0}^{(p-3)/2}\f{\bi{2k}k}{(2k+1)16^k}\eq0\ (\mo\ p^2)\  \t{and} \ \sum_{k=(p+1)/2}^{p-1}\f{\bi{2k}k}{(2k+1)16^k}\eq\f p3E_{p-3}\ (\mo\ p^2),$$
and conjectured that
$$\sum_{k=0}^{(p-3)/2}\f{\bi{2k}k}{(2k+1)^2(-16)^k}\eq\f{H_{p-1}}{5p}\pmod{p^3}$$
(which was recently confirmed by K. Hessami Pilehrood, T. Hessami Pilehrood and R. Tauraso [HPT]) and
$$\sum_{k=(p+1)/2}^{p-1}\f{\bi{2k}k}{(2k+1)^2(-16)^k}\eq-\f p4B_{p-3}\pmod{p^2},$$
where $p$ is any prime greater than 5.

Theorem 1.1 has the following consequence.

\proclaim{Corollary 1.1} Let $p>5$ be a prime. Then
$$\f1p\sum_{k=(p+1)/2}^{p-1}\f{\bi{2k}k}{(2k+1)^2(-16)^k}\eq-\sum_{k=0}^{(p-3)/2}\f{(-16)^k}{(2k+1)^3\bi{2k}k}\eq-\f {B_{p-3}}4\pmod{p}.\tag1.5$$
\endproclaim

Since
$$\sum_{k=0}^n\f{\bi{2k}k^2}{(2k-1)16^k}=-\f{2n+1}{16^n}\bi{2n}n^2\ \t{and}\ \sum_{k=0}^n\f{(1-4k)\bi{2k}k^4}{(2k-1)^4256^k}=(8n^2+4n+1)\f{\bi{2n}n^4}{256^n}$$
by induction, we find that
$$\sum_{k=0}^\infty\f{\bi{2k}k^2}{(2k-1)16^k}=-\f 2{\pi}\ \ \t{and}\ \ \sum_{k=0}^\infty\f{(4k-1)\bi{2k}k^4}{(2k-1)^4256^k}=-\f 8{\pi^2}$$
by Stirling's formula $n!\sim\sqrt{2\pi n}(n/e)^n$. (The latter was first obtained by J. W. L. Glaisher [G2].)
Via {\tt Mathematica} (version 7) we find that
$$\sum_{k=0}^\infty\f{\bi{2k}k^2}{(2k+1)16^k}=\f{4G}{\pi}
\ \ \t{and}\ \ \sum_{k=0}^\infty\f{16^k}{(2k+1)^3\bi{2k}k^2}=\f 72\zeta(3)-G\pi.$$
Motivated by this we establish the following theorem.

\proclaim{Theorem 1.2} Let $p>3$ be a prime. Then
$$\f1{p^2}\sum_{k=(p+1)/2}^{p-1}\f{\bi{2k}k^2}{(2k+1)16^k}\eq-\sum_{k=0}^{(p-3)/2}\f{16^k}{(2k+1)^3\bi{2k}k^2}\eq-\f 74B_{p-3}\pmod p,\tag1.6$$
and
$$\sum_{k=0}^{(p-3)/2}\f{\bi{2k}k^2}{(2k+1)16^k}\eq-2q_p(2)-p\,q_p(2)^2+\f 5{12}p^2B_{p-3}\pmod{p^3},\tag1.7$$
where $q_p(2)$ denotes the Fermat quotient $(2^{p-1}-1)/p$.
\endproclaim

Now we pose two conjectures for further research.

\proclaim{Conjecture 1.1} {\rm (i)}
If $p>5$ is a prime, then
$$\align\sum_{k=(p+1)/2}^{p-1}\f{\bi{2k}k^2}{k16^k}&\eq-\f {21}2H_{p-1}\ (\mo\ p^4),
\\\sum_{k=0}^{(p-3)/2}\f{(-16)^k}{(2k+1)^3\bi{2k}k}&\eq-\f{3}4\cdot\f{H_{p-1}}{p^2}-\f{47}{400}p^2B_{p-5}\pmod{p^3}.
\endalign$$

{\rm (ii)} If $p>3$ is a prime, then
$$\align\sum_{k=1}^{(p-1)/2}\f{\bi{2k}k^2H_{2k}}{k16^k}&\eq (-1)^{(p-1)/2}4E_{p-3}\pmod{p},
\\\sum_{k=1}^{(p-1)/2}\f{\bi{2k}k^2}{k16^k}(H_{2k}-H_k)&\eq-\f73pB_{p-3}\pmod{p^2},
\\p^2\sum_{k=1}^{p-1}\f{16^{k-1}H_{k-1}}{k^2\bi{2k}k^2}&\eq\f{(-1)^{(p-1)/2}}2H_{(p-1)/2}+pE_{p-3}\pmod{p^2}.
\endalign$$ We also have the identity
$$\sum_{k=1}^\infty\f{\bi{2k}k^2}{k16^k}(H_{2k}-H_k)=\f23\sum_{k=0}^\infty\f{\bi{2k}k^2H_{2k}}{(2k+1)16^k}.\tag1.8$$
\endproclaim

\proclaim{Conjecture 1.2} Let $p$ be an odd prime. Then
$$\align\sum_{k=1}^{p-1}\f{\bi{2k}k\bi{4k}{2k}}{k64^k}\eq&-3H_{(p-1)/2}+\f 74p^2B_{p-3}\ (\mo\ p^3),
\\\sum_{k=1}^{(p-1)/2}\f{\bi{2k}k\bi{4k}{2k}}{k64^k}\eq&-3H_{(p-1)/2}+(-1)^{(p+1)/2}\,2pE_{p-3}\ (\mo\ p^2),
\\p\sum_{k=1}^{(p-1)/2}\f{64^{k-1}}{k^3\bi{2k}k\bi{4k}{2k}}\eq&\f{(-1)^{(p-1)/2}}{2}E_{p-3}\ \ (\mo\ p).
\endalign$$
If $p>3$, then
$$p\sum_{k=1}^{(p-1)/2}\f{64^{k-1}}{(2k-1)k^2\bi{2k}k\bi{4k}{2k}}\eq \f{(-1)^{(p+1)/2}q_p(2)+pE_{p-3}}4\ \ (\mo\ p^2)$$
and
$$\sum_{k=1}^{(p-1)/2}\f{\bi{2k}k\bi{4k}{2k}}{k64^k}(H_{2k}-H_k)\eq(-1)^{(p+1)/2}4E_{p-3}\pmod{p}.$$
\endproclaim

\medskip
In the next section we are going to show Theorem 1.1 and Corollary
1.1. Theorem 1.2 will be proved in Section 3. Our proofs involve
certain combinatorial identities and harmonic numbers of higher
orders given by
$$H_n^{(m)}=\sum_{0<k\ls n}\f1{k^m}\ \ \ (n=0,1,2,\ldots).$$

\heading{2. Proof of Theorem 1.1}\endheading

\proclaim{Lemma 2.1} Let $p=2n+1$ be an odd prime. For any $k=1,\ldots,p-1$, we have
$$k\bi{2k}k\bi{2(p-k)}{p-k}\eq (-1)^{\lfloor 2k/p\rfloor-1}2p\pmod{p^2}\tag2.1$$
and
$$\bi nk\bi{n+k}k\eq\f{\bi{2k}k^2}{(-16)^k}\pmod{p^2}.\tag2.2$$
\endproclaim
\Proof. Congruence (2.1) was formulated in [S11c, Lemma 2.1]; see
also [T10] for such a trick. Congruence (2.2) is known and easy (see
p.\,231 of [vH, \S3.4]); in fact,
$$\bi nk\bi{n+k}k(-1)^k=\bi{(p-1)/2}k\bi{(-p-1)/2}k\eq\bi{-1/2}k^2=\f{\bi{2k}k^2}{16^k}\ (\mo\ p^2).$$
We are done. \qed

\proclaim{Lemma 2.2} Let $p>3$ be a prime. Then
$$(-1)^{(p-1)/2}\bi{p-1}{(p-1)/2}\eq 4^{p-1}+\f{p^3}{12}B_{p-3}\pmod{p^4}\tag2.3$$
and
$$H_{(p-1)/2}\eq-2q_p(2)+p\,q_p(2)^2-p^2\l(\f23q_p(2)^3+\f 7{12}B_{p-3}\r)\pmod{p^3}.\tag2.4$$
Also,
$$H_{(p-1)/2}^{(2)}\eq\f 73pB_{p-3}\ (\mo\ p^2)\ \ \t{and}\ \ H_{(p-1)/2}^{(3)}\eq-2B_{p-3}\ (\mo\ p).\tag2.5$$
\endproclaim
\Remark\ 2.1. (2.3) and (2.4) are refinements of Morley's congruence [Mo] and Lehmer's congruence [L] given by L. Carlitz [C]
and Z.-H. Sun [Su, Theorem 5.2(c)] respectively. (2.5) follows from [Su, Corollary 5.2].

\proclaim{Lemma 2.3} For each $n=1,2,3,\ldots$ we have
$$\sum_{k=1}^n\f{(-1)^k}{k^3\bi nk\bi{n+k}k}=5\sum_{k=1}^n\f{(-1)^k}{k^3\bi{2k}k}+2H_n^{(3)}\tag2.6$$
and
$$\sum_{k=1}^n\bi nk\bi{n+k}k\f{(-1)^k}k(H_{n+k}-H_{n-k})=\f 52\sum_{k=1}^n\f{(-1)^k}{k^2}\bi{2k}k+2H_n^{(2)}.\tag2.7$$
\endproclaim
\Proof. (2.6) is due to Ap\'ery [Ap] (see also [vP]). The author
found (2.7) via the math. software {\tt Sigma}. (The reader may
consult [OS, \S5] for how to use {\tt Sigma} to produce
combinatorial identities.) In fact, if we let $s_n$ denote the
left-hand side or the right-hand side of (2.7) then {\tt Sigma}
yields the recurrence relation
$$(n+1)^2(s_{n+1}-s_n)=2-5(-1)^n\bi{2n+1}n\ \ (n=1,2,3,\ldots).$$
So (2.7) can be proved by induction. \qed

\proclaim{Lemma 2.4} Let $p>3$ be a prime. Then
$$\sum_{k=1}^{p-1}\f{\bi{2k}k^2}{k^216^k}\eq -2H_{(p-1)/2}^2\pmod{p^2}$$
and
$$\sum_{k=1}^{p-1}\f{\bi{2k}k^2}{k^3 16^k}\eq -\f 43H_{(p-1)/2}^3-\f 23H_{(p-1)/2}^{(3)}\pmod{p}.$$
\endproclaim
\Remark\ 2.2. Lemma 2.4 is known, see [T12, Theorem 7] and its proof.

\proclaim{Lemma 2.5} We have the new identity
$$\sum_{k=0}^n\f{\bi{2k}k^2}{(2(n+k)+1)16^k}=\f{\bi{2n}n^2}{16^n}\sum_{k=0}^{2n}\f1{2k+1}.\tag2.8$$
\endproclaim
\Proof. Let $u_n$ denote the left-hand side or the right-hand side of (2.8).
Applying the Zeilberger algorithm (cf. [PWZ, pp.\,101-119]) via {\tt Mathematica} (version 7), we find the recurrence relation
$$(2n+1)^2u_n-4(n+1)^2u_{n+1}=-\f{8(4n^3+8n^2+5n+1)}{(4n+3)(4n+5)16^n}\bi{2n}n^2\ \ (n=0,1,2,\ldots).$$
 So (2.8) holds by induction. \qed

\medskip
 \noindent{\it Proof of Theorem 1.1}. For convenience we set $n=(p-1)/2$ and $S:=\sum_{k=1}^n(-1)^k/(k^3\bi{2k}k)$.
 Below we divide the proof into two parts.
\medskip

 (i) In light of Lemma 2.1,
 $$\f1p\sum_{k=n+1}^{p-1}\f{(-1)^k}{k^2}\bi{2k}k\eq\sum_{k=n+1}^{p-1}\f{2(-1)^k}{k^3\bi{2(p-k)}{p-k}}
 =\sum_{k=1}^n\f{2(-1)^{p-k}}{(p-k)^3\bi{2k}k}\eq2S\pmod p.$$
Hence (1.1) and (1.2) are equivalent since
$$\sum_{k=1}^{p-1}\f{(-1)^k}{k^2}\bi{2k}k\eq -\f 4{15}pB_{p-3}\pmod{p^2}\tag2.9$$
by [T10].

By [S11c, (3.2)], for $k=1,\ldots,n$ we have
$$\bi nk\bi{n+k}k(-1)^k\l(1-\f p4(H_{n+k}-H_{n-k})\r)\eq\f{\bi{2k}k^2}{16^k}\pmod{p^4}.$$
This, together with (2.7) and the known identity
$$\sum_{k=1}^n\bi nk\bi{n+k}k\f{(-1)^k}k=-2H_n$$
(cf. [Pr, \S2.1]), yields that
$$\sum_{k=1}^n\f{\bi{2k}k^2}{k16^k}+2H_n
\eq-\f p4\(\f 52\sum_{k=1}^n\f{(-1)^k}{k^2}\bi{2k}k+2H_n^{(2)}\)\pmod{p^4}.$$
Thus, in view of the first congruence in (2.5), we have
$$\sum_{k=1}^n\f{\bi{2k}k^2}{k16^k}+2H_n\eq-\f 58p\sum_{k=1}^n\f{(-1)^k}{k^2}\bi{2k}k-\f 76p^2B_{p-3}\pmod{p^3}.\tag2.10$$
Therefore (1.2) and (1.3) are equivalent.

In light of Lemma 2.1,
 $$\align\f1{p^2}\sum_{k=n+1}^{p-1}\f{\bi{2k}k^2}{k16^k}\eq&\sum_{k=n+1}^{p-1}\f 4{k^316^k\bi{2(p-k)}{p-k}^2}
 =\sum_{k=1}^n\f 4{(p-k)^3 16^{p-k}\bi{2k}k^2}
 \\\eq&-\f14\sum_{k=1}^n\f{16^k}{k^3\bi{2k}k^2}\pmod{p}.
 \endalign$$
With the help of (2.2) and (2.5), we obtain from (2.6) the congruence
$$\sum_{k=1}^n\f{16^k}{k^3\bi{2k}k^2}\eq5S-4B_{p-3}\pmod p.$$
Therefore (1.1) and (1.4) are equivalent.

Note that
$$\sum_{k=n+1}^{p-1}\f{\bi{2k}k^2}{k16^k}\eq-\f{p^2}4(5S-4B_{p-3})\eq-\f 58p\sum_{k=n+1}^{p-1}\f{(-1)^k}{k^2}\bi{2k}k+p^2B_{p-3}\pmod{p^3}.$$
Combining this with (2.10) we get
$$\sum_{k=1}^{p-1}\f{\bi{2k}k^2}{k16^k}+2H_n\eq-\f 58p\sum_{k=1}^{p-1}\f{(-1)^k}{k^2}\bi{2k}k-\f{p^2}6B_{p-3}\pmod{p^3}.\tag2.11$$
So Tauraso's result $\sum_{k=1}^{p-1}\bi{2k}k^2/(k16^k)\eq-2H_n\ (\mo\ p^3)$ in [T12] actually follows from his earlier result (2.9).
\medskip

(ii) By the above, it suffices to show (1.3). Note that
$$\sum_{k=0}^{2n}\f1{2k+1}=H_{4n+1}-\f{H_{2n}}2=H_{p-1}+\f1{p}+\sum_{k=1}^{p-1}\f1{p+k}-\f{H_{p-1}}2$$
and
$$\align\sum_{k=1}^{p-1}\f1{p+k}=&\sum_{k=1}^{p-1}\f{k^2-kp+p^2}{k^3+p^3}
\\\eq&\sum_{k=1}^{p-1}\f{k^2-kp+p^2}{k^3}=H_{p-1}-pH_{p-1}^{(2)}+p^2H_{p-1}^{(3)}\pmod{p^3}.\endalign$$
As
$$H_{p-1}^{(3)}=\sum_{k=1}^{(p-1)/2}\l(\f1{k^3}+\f1{(p-k)^3}\r)\eq0\pmod p,$$
we see that
$$\sum_{k=0}^{2n}\f1{2k+1}-\f1p\eq\f 32H_{p-1}-pH_{p-1}^{(2)}\eq-\f {p^2}2B_{p-3}-\f 23p^2B_{p-3}=-\f 76p^2B_{p-3}\ (\mo\ p^3)$$
since $H_{p-1}\eq-p^2B_{p-3}/3\ (\mo\ p^3)$ and $H_{p-1}^{(2)}\eq2pB_{p-3}/3\ (\mo\ p^2)$ by Glaisher [G1] (see also [Su, Theorem 5.1(a)]).
In view of (2.3),
$$\align\f{\bi{2n}n^2}{16^n}\eq&\f{(4^{p-1}(1+p^3B_{p-3}/12))^2}{4^{p-1}}\eq 4^{p-1}\l(1+\f{p^3}6B_{p-3}\r)
\\\eq& (1+p\,q_p(2))^2+\f{p^3}6B_{p-3}\pmod{p^4}.
\endalign$$
Therefore
$$\align\f{\bi{2n}n^2}{16^n}\sum_{k=0}^{2n}\f1{2k+1}-\f1p=&\f{\bi{2n}n^2}{16^n}\(\sum_{k=0}^{2n}\f1{2k+1}-\f1p\)+\f{\bi{2n}n^2/16^n-1}p
\\\eq&-\f 76p^2B_{p-3}+\f{(1+p\,q_p(2))^2+p^3B_{p-3}/6-1}p
\\=& 2q_p(2)+p\,q_p(2)^2-p^2B_{p-3}\pmod{p^3}
\endalign $$
and hence (2.8) yields the congruence
$$\sum_{k=1}^n\f{\bi{2k}k^2}{(2k+p)16^k}\eq 2q_p(2)+p\,q_p(2)^2-p^2B_{p-3}\pmod{p^3}.\tag2.12$$

Note  that
$$\aligned\sum_{k=1}^n\f{\bi{2k}k^2}{(2k+p)16^k}=&\sum_{k=1}^n\f{(4k^2-2kp+p^2)\bi{2k}k^2}{((2k)^3+p^3)16^k}
\\\eq&\f12\sum_{k=1}^n\f{\bi{2k}k^2}{k16^k}-\f p4\sum_{k=1}^n\f{\bi{2k}k^2}{k^216^k}+\f{p^2}8\sum_{k=1}^n\f{\bi{2k}k^2}{k^316^k}\pmod{p^3}.
\endaligned\tag2.13$$
By Lemma 2.4 and (2.4)-(2.5),
$$\sum_{k=1}^n\f{\bi{2k}k^2}{k^216^k}\eq-2(-2q_p(2)+p\,q_p(2)^2)^2\eq-8q_p(2)^2+8p\,q_p(2)^3\pmod{p^2}\tag2.14$$
and
$$\sum_{k=1}^n\f{\bi{2k}k^2}{k^316^k}\eq-\f 43(-2q_p(2))^3-\f 23(-2B_{p-3})=\f{32}3q_p(2)^3+\f 43B_{p-3}\pmod{p}.\tag2.15$$
Combining (2.12)-(2.15) we obtain
$$\align\f 12\sum_{k=1}^n\f{\bi{2k}k^2}{k16^k}
\eq&2q_p(2)+p\,q_p(2)^2-p^2B_{p-3}+\f p4(-8q_p(2)^2+8p\,q_p(2)^3)
\\&-\f{p^2}8\l( \f{32}3q_p(2)^3+\f 43B_{p-3}\r)\pmod{p^3}.
\endalign$$
In view of (2.4), this yields the desired (1.3).

The proof of Theorem 1.1 is now complete. \qed

\medskip
\noindent{\it Proof of Corollary 1.1}. Write $p=2n+1$. In view of Lemma 2.1,
$$\align &\f1p\sum_{k=n+1}^{p-1}\f{\bi{2k}k}{(2k+1)^2(-16)^k}
\\=&\sum_{k=1}^n\f{\bi{2(p-k)}{p-k}/p}{(2(p-k)+1)^2(-16)^{p-k}}
\\\eq&\sum_{k=1}^n\f{-2/(k\bi{2k}k)}{(2(k-1)+1)^2(-16)^{1-k}}=-\sum_{k=0}^{n-1}\f{(-16)^k}{(2k+1)^3\bi{2k}k}\pmod{p}.
\endalign$$
Since
$$\bi{-1/2}k\eq\bi nk=\bi n{n-k}\eq\bi{-1/2}{n-k}\pmod p\quad\t{for all}\ k=0,\ldots,n,$$
we have
$$\align\sum_{k=0}^{n-1}\f{(-16)^k}{(2k+1)^3\bi{2k}k}=&\sum_{k=0}^{n-1}\f{4^k}{(2k+1)^3\bi{-1/2}k}
\\\eq&\sum_{k=0}^{n-1}\f{4^k}{(2k+1)^3\bi{-1/2}{n-k}}=\sum_{k=1}^n\f{4^{n-k}}{(2(n-k)+1)^3\bi{-1/2}k}
\\\eq&-\f18\sum_{k=1}^n\f1{k^3\bi{-1/2}k4^k}=-\f18\sum_{k=1}^n\f{(-1)^k}{k^3\bi{2k}k}
\\\eq&\f{B_{p-3}}4\pmod{p}\quad\t{(by (1.1))}.
\endalign$$
Therefore (1.5) holds. \qed

\heading{3. Proof of Theorem 1.2}\endheading

\proclaim{Lemma 3.1} For any positive integer $n$ we have
$$\sum_{k=0}^{n-1}\f{\bi{2k}k^2}{(n-k)16^k}=\f{\bi{2n}n^2}{4^{2n-1}}\sum_{k=0}^{n-1}\f1{2k+1}\tag3.1$$
and
$$\sum_{k=0}^{n}\f{(-1)^k}{(2k+1)^2}\bi nk\bi{n+k}k=\f{1}{(2n+1)^2}+\f2{2n+1}\sum_{k=0}^{n-1}\f1{2k+1}.\tag3.2$$
\endproclaim
\Proof. Let $(x)_0=1$ and $(x)_k=\prod_{j=0}^{k-1}(x+j)$ for $k=1,2,3,\ldots$. Then
$(1/2)_k^2/(1)_k^2=\bi{2k}k^2/16^k$ for $k=0,1,2,\ldots$. Thus (3.1) is just (21) of [Lu, Ch.\,5.2] with $x=1/2$
(see also [T12, (1)]).

Identity (3.2) is new and it can be proved via the Zeilberger algorithm (cf. [PWZ, pp.\,101--119]). It is easy to verify (3.2) for $n=1,2$.
Applying the Zeilberger algorithm via {\tt Mathematica} (version 7), we find that if $a_n$ denotes the left-hand side or the right-hand side of (3.2)
then
$$(n+1)(2n+5)^2a_{n+2}=(2n+3)(4n^2+12n+7)a_{n+1}-(n+2)(2n+1)^2a_n$$
for all $n=1,2,3,\ldots$. So (3.2) follows by induction. \qed

\proclaim{Lemma 3.2} Let $p>3$ be a prime. Then
$$\sum_{k=0}^{(p-3)/2}\f{\bi{2k}k^2}{(2k+1)^316^k}\eq-\f 43q_p(2)^3-\f{B_{p-3}}6\pmod p.\tag3.3$$
\endproclaim
\Proof. Set $n=(p-1)/2$. For $k=0,\ldots,n$ we clearly have
$$\bi{n+k}k(-1)^k=\bi{-n-1}k\eq\bi{n}k\eq\bi{-1/2}k=\f{\bi{2k}k}{(-4)^k}\pmod p.$$
Thus
$$\align\sum_{k=0}^{n-1}\f{\bi{2k}k^2}{(2k+1)^316^k}
\eq&\sum_{k=0}^{n-1}\f{\bi nk^2}{(2k+1)^3}=\sum_{k=1}^n\f{\bi nk^2}{(2(n-k)+1)^3}
\\\eq&-\f18\sum_{k=1}^{n}\f{\bi{2k}k^2}{k^316^k}\pmod p.
\endalign$$
Combining this with (2.15) we obtain the desired (3.3).
\qed

\proclaim{Lemma 3.3} Let $p>3$ be a prime. Then
$$\sum_{k=0}^{(p-3)/2}\f{\bi{2k}k^2}{(2k+1)^216^k}\eq-2q_p(2)^2+\f 23p\,q_p(2)^3-\f p6B_{p-3}\pmod{p^2}.\tag3.4$$
\endproclaim
\Proof. Write $p=2n+1$. By (3.2) we have
$$\align\sum_{k=0}^{n-1}\f{(-1)^k}{(2k+1)^2}\bi nk\bi{n+k}k=&\f{1-(-1)^n\bi{2n}n}{p^2}+\f2{p}\l(H_{p-1}-\f{H_n}2\r)
\\=&\f{1-(-1)^n\bi{2n}n-pH_n}{p^2}+\f 2pH_{p-1}.
\endalign$$
In light of Lemma 2.2,
$$\align1-(-1)^n\bi{2n}n-pH_n
\eq&1-(1+p\,q_p(2))^2-\f{p^3}{12}B_{p-3}
\\&+2p\,q_p(2)-p^2q_p(2)^2+p^3\l(\f 23q_p(2)^3+\f7{12}B_{p-3}\r)
\\=&-2p^2q_p(2)^2+p^3\l(\f 23q_p(2)^3+\f{B_{p-3}}2\r)\pmod{p^4}.
\endalign$$
So, with the help of (2.2), we get
$$\sum_{k=0}^{n-1}\f{\bi{2k}k^2}{(2k+1)^2 16^k}\eq-2q_p(2)^2+p\l(\f 23q_p(2)^3+\f{B_{p-3}}2\r)+\f2pH_{p-1}\pmod{p^2},$$
which gives (3.4) since $H_{p-1}/p\eq -pB_{p-3}/3\ (\mo\ p^2)$. \qed

\medskip
\noindent{\it Proof of Theorem 1.2}. For convenience we set $n=(p-1)/2$.

We first prove (1.6). By Lemma 2.1,
$$\align&\f1{p^2}\sum_{k=n+1}^{p-1}\f{\bi{2k}k^2}{(2k+1)16^k}
=\sum_{k=1}^n\f{(\bi{2(p-k)}{p-k}/p)^2}{(2(p-k)+1)16^{p-k}}
\\\eq&\sum_{k=1}^n\f{(-2/(k\bi{2k}k))^216^{k-1}}{1-2k}=-\sum_{k=0}^{n-1}\f{16^k}{(2k+1)^3\bi{2k}k^2}\pmod p.
\endalign$$
Note also that
$$\align \sum_{k=0}^{n-1}\f{16^k}{(2k+1)^3\bi{2k}k^2}\eq&\sum_{k=0}^{n-1}\f1{(2k+1)^3\bi nk^2}=\sum_{k=1}^n\f1{(2(n-k)+1)^3\bi nk^2}
\\\eq&-\f18\sum_{k=1}^n\f{1}{k^3\bi nk^2}\eq-\f18\sum_{k=1}^n\f{16^k}{k^3\bi{2k}k^2}\pmod p.
\endalign$$
Thus, with the help of (1.4) we get (1.6). (In the case $p=5$, (1.6) can be verified directly.)

Since $\bi{2n}n^2\eq 4^{4n}\pmod{p^3}$ by (2.3), we have
$$\align\f{\bi{2n}n^2}{4^{2n-1}}\sum_{k=0}^{n-1}\f1{2k+1}\eq &4^p\l(H_{p-1}-\f{H_n}2\r)=4^pH_{p-1}-2(1+p\,q_p(2))^2H_n
\\\eq&-\f {4}3p^2B_{p-3}-2(1+2p\,q_p(2)+p^2q_p(2)^2)H_n
\\\eq&4q_p(2)+6p\,q_p(2)^2+p^2\l(\f 43q_p(2)^3-\f{B_{p-3}}6\r)\pmod{p^3}
\endalign$$
with the help of (2.4). Note also that
$$\align &\f12\sum_{k=0}^{n-1}\f{\bi{2k}k^2}{(n-k)16^k}=\sum_{k=0}^{n-1}\f{\bi{2k}k^2}{(p-(2k+1))16^k}
\\=&\sum_{k=0}^{n-1}\f{\bi{2k}k^2(p^2+p(2k+1)+(2k+1)^2)}{(p^3-(2k+1)^3)16^k}
\\\eq&-p^2\sum_{k=0}^{n-1}\f{\bi{2k}k^2}{(2k+1)^316^k}-p\sum_{k=0}^{n-1}\f{\bi{2k}k^2}{(2k+1)^216^k}
-\sum_{k=0}^{n-1}\f{\bi{2k}k^2}{(2k+1)16^k}\pmod{p^3}.
\endalign$$
Combining these with (3.1), (3.3) and (3.4) we finally obtain the desired (1.7). \qed

\Ack. The author is grateful to the referee for helpful comments.

\enddocument